 \documentclass[12pt]{article}
 \usepackage{amsmath,amsfonts,theorem}
 \numberwithin{equation}{section}

 \newtheorem{prop}{Proposition}[section]

 \newtheorem{cor}{Corollary}[section]

 \newtheorem{thm}{Theorem}[section]

 \theorembodyfont{\rmfamily}
 \newtheorem{dfn}{Definition}[section]

 \newtheorem{rmk}{Remark}

 \newcommand{\qed}{\ifhmode\unskip\nobreak\fi\quad\ensuremath\square}

 \newcommand{\CS}{\operatorname{CS}}

 \newcommand{\semis}{^{\mathrm{ss}}}

 \newcommand{\rest}[1]{{}_{{\textstyle{|}}#1}} 

 \newcommand{\p}{\partial}

 \newcommand{\ov}{\overline}

 \newcommand{\sA}{\mathcal A} 

 \newcommand{\sG}{\mathcal G} 

 \newcommand{\sH}{\mathcal H}
 \newcommand{\sL}{\mathcal L}
 \newcommand{\sM}{\mathcal M}
 \newcommand{\Oh}{\mathcal O}

 \newcommand{\al}{\alpha}

 \newcommand{\om}{\omega}

 \newcommand{\Ga}{\Gamma}
 \newcommand{\La}{\Lambda}
 \newcommand{\Om}{\Omega}
 
 \newcommand{\Si}{\Sigma}
 \newcommand{\na}{\nabla}

 \newcommand{\PP}{\mathbb P}
 \newcommand{\C}{\mathbb C}

 \newcommand{\R}{\mathbb R}
 \newcommand{\Z}{\mathbb Z}

 \newcommand{\rk}{\operatorname{rank}}

 \newcommand{\Diff}{\operatorname{Diff}}

 \newcommand{\tr}{\operatorname{tr}}

 \newcommand{\Hom}{\operatorname{Hom}}
 
 \newcommand{\PU}{\operatorname{PU}}
 \newcommand{\SU}{\operatorname{SU}}

\begin{document}

 \title{Complexification  of  Bohr-Sommerfeld conditions.}
 \markright{\hfill Complexification  of BS-conditions \quad}

 \author{Andrei Tyurin}
 \date{ September,13, 1999}
 \maketitle

 \begin{abstract}
The complex version of Bohr-Sommerfeld conditions is proposed. The BPU-construction  (see \cite{BPU} or \cite{T1}) is generalized to this complexification. The new feature of this generalization is a spectral curve. The geometry of such curves is
 investigated.  
\end{abstract}

 \section{Global structures on  subspaces of Lagrangian cycles}

Let $(M, \om)$ be a smooth symplectic manifold of dimesion $2n$. It can be considered as a phase space of a classical mechanical system. The lagrangian submanifolds play an especially important role in symplectic geometry. Let $\sL$ be a smooth manifold of dimension $n$ and 
\begin{equation}
i \colon \sL \to M
\label{eq1.1}
\end{equation}
be  an embedding of $\sL$ into $M$. Then  we can consider the subcycle $i(\sL) \subset M $  as an orbit in the space $Map (\sL, M)$ of smooth maps to $M$ with respect to the natural action of the diffeomorphisms group $\Diff (\sL)$:
\begin{equation}
(\sL \subset M) = Map (\sL, M) / \Diff (\sL).
\label{eq1.2}
\end{equation}
Such a  submanifold is called Lagrangian if 
\begin{equation}
i^* (\om) = 0.
\label{eq1.3}
\end{equation}
It is easy to see that the normal bundle
\begin{equation}
N_M \sL = T^* \sL
\label{eq1.4}
\end{equation}
-the cotangent bundle of $\sL$.
By the Darboux-Weinstein theorem a Lagrangian submanifold $(\sL \subset M)$ has not invariants of embeddings: a small tubular neighborhood of $ \sL$ can be identified with the  neighborhood of the zero section of the cotangent bundle. Moreover, any Lagrangian subcycle in this  neighborhood can be identified with a closed 1-form on $\sL$. Hence such closed forms define a chart around a submanifold $(\sL \subset M)$ in the space $\sL(M)$ of all Lagrangian subcycles of $M$. Moreover
 the tangent space at $\sL$
\begin{equation}
T \sL(M)_\sL = \{ \al \in \Om(\sL) \quad \vert \quad d \al = 0 \} = Z^1(\sL)
\label{eq1.5}
\end{equation}
- the space of closed 1-forms on $\sL$. 

So besides the  smooth type of $\sL$,  Lagrangian submanifolds look like points of $M$.

On the other hand we can consider the trivial Hermitian line bundle $ L_0 $
 on $\sL$. Then every closed form $\al \in Z^1(\sL)$ multiplied by $ i = \sqrt{-1}$  can be  consided  as a flat connection  on $L_0$ and the gauge class of this connection defines  a character of the fundamental group
\begin{equation}
\chi_\al \colon \pi_1 (\sL) \to U(1) \in H^1(\sL, \R) / H^1(\sL, \Z) = J_\sL.
\label{eq1.6}
\end{equation} 
(Let us call the last torus  the Jacobian of $\sL$).

Thus   an infinitesimal deformation $\al \in Z^1$ of $\sL$ defines the point
\begin{equation}
 [\al ] \in J_\sL
\label{eq1.7}
\end{equation}
and the space of infinitesimal deformation of $\sL$  preserving  some fixed class (for example, trivial) of  flat connections
 is
\begin{equation}
 B^1(\sL) = \{ \al \in Z^1(\sL) \quad \vert \quad \al = d f \}
\label{eq1.8}
\end{equation}
where $ f \in C^\infty (\sL)$. Such deformations of Lagrangian submanifolds are called {\it isodrastic } deformations (see \cite{Wei2} for the linguistic  explanation of this term). So the tangent space to isodrastic deformations 
\begin{equation}
 T_I \sL(M)_\sL = B^1(\sL) = C^\infty (\sL) / \R.
\label{eq1.9}
\end{equation}

 To get some interesting structures on the space of all Lagrangian cycles $\sL (M)$ we have to equip our Lagrangian submanifolds with some extra structure. We will consider connected Lagrangian cycles only.

Recall that a {\it half-weithed } Lagrangian submanifold  is a pair $(\sL, hF)$ where $hF$ is a smooth half-form on $\sL$. A half-weighted manifold is {\it weighted} automatically because $hF^2 = V$ is a volume form on $\sL$. Hence we have the double cover 
\begin{equation}
sq \colon hW\sL(M) \to W\sL(M)
\label{eq1.10}
\end{equation}
of the space of {\it  weighted} Lagrangian manifolds ramified along 
\[
\sL(M) = \{ (\sL, hF) \quad \vert \quad hF = 0 \}.
\]
 Thus the tangent space at $( \sL_0, hF_0 )$
\begin{equation}
T hW\sL(M)_{( \sL_0, hF_0 )} = Z^1(\sL) \oplus  C^\infty (\sL) \cdot hF_0
\label{eq1.11}
\end{equation}
where the last vector space we consider as a space of half-forms. This space carries  the form
\begin{equation}
(f_1 \cdot hF_0, f_2 \cdot hF_0) = \int_\sL f_1 \cdot f_2 \cdot hF^2.
\label{eq1.12}
\end{equation}
 
The space
$hW\sL(M)$ of all half-weighted Lagrangian submanifolds of $M$ 
admits the function 
\begin{equation}
v \colon hW\sL(M)  \to        \R , \quad v (\sL, hF) = \int_\sL hF^2.
\label{eq1.13}
\end{equation}
Consider  level hypersurfaces of this function 
\begin{equation}
 hW\sL(M)_t = v^{-1}(t)
\label{eq1.14}
\end{equation}
-the space of half-weighted Lagrangian cycles of  volume $t$.

 Thus when we deform $( \sL_0, hF_0 )$ inside of  $ hW\sL(M)_t$ then by Moser theorem \cite{Wei1} (more precisely, by its  paraphrase ) each half-weighted submanifold $(\sL, hF)$ in the family of deformation is diffeomorphic as a half-weighted manifold with  $( \sL_0, hF_0 )$. For example
\begin{equation}
 (T hW\sL(M)_0)_{( \sL_0, hF_0 )} = Z^1(\sL) \oplus (hF_0)^\perp
\label{eq1.15}
\end{equation}
-the orthogonal subspace with respect to the form (1.8).

On the  subspace of isodrastic deformations of a  half-weighted submanifold 
\begin{equation}
(T_I hW\sL(M)_0)_{( \sL_0, hF_0 )} = B^1(\sL) \oplus (hF_0)^\perp
\label{eq1.16}
\end{equation}
we can define  a 2-form by 
\begin{equation}
\Om_0 ((f_1, hF_1), (f_2, hF_2)) = \int_\sL (f_1 \cdot hF_2 - f_2 \cdot hF_1) \cdot hF_0. 
\label{eq1.17}
\end{equation}
It's easy to see that this form is weak non degenerated. Moreover:

\begin{prop} This form is closed.
\end{prop}

We can thus think of the space of isodrastic deformations of half-weighted Lagrangian  submanifolds of fixed volume  as an infinite-dimensional symplectic manifold.

In the same vein we get the symplectic form on the space
\begin{equation}
IhW\sL(M)_t \subset IhW\sL(M)
\label{eq1.18}
\end{equation}
of isodrastic deformations of  half-weighted  cycles of volume t (see below).

The next extra structure is coming from

\subsection*{Prequantization}

Let us equip our symplectic manifold $M$ with a complex line bundle $L$  carring a Hermitian connection $a_L$ such that its curvature form
\begin{equation}
F_{a_L} = 2 \pi i \cdot k \cdot \om, \quad k \in \Z.
\label{eq1.19}
\end{equation}
Such line bundle $(L, a_L)$ is called a {\it prequantization} bundle of the phase space $(M, \om)$ of  a mechanical system and the integer $k$ is called a {\it level}. Such quadruple
\[
(M, \om, L, a_L)
\]
is a phase space of a mechanical system which is ready to be quantized.

For every Lagrangian submanifold $\sL \in \sL(M)$ the  restriction of $ (L, a_L)$ is the trivial line bundle with a  flat connection. A gauge class of this connection defines a point
\begin{equation}
\al( a_L) \in J_\sL
\label{eq1.20}
\end{equation}
on the Jacobian of $\sL$ ( see (1.7)). So we can distinguish Lagrangian submanifolds by these points.

\begin{dfn} A Lagrangian cycle $\sL$ is called {\it Bohr-Sommerfeld} (BS for short) if
\[
\al( a_L) = 0  
\] 
(see (1.7), (1.8) and (1.9)).
\end{dfn}
We get the submanifold
\begin{equation}
BS( a_L) \subset \sL(M)  
\label{eq1.21}
\end{equation} 
of Bohr-Sommerfeld with respect to a  pair $(L, a_L)$. 

This subspace is isodrastic. The tangent space at $\sL \in BS( a_L)$ is
\begin{equation}
T BS( a_L)_\sL = B^1(\sL)
\label{eq1.22}
\end{equation}
( see (1.8) and (1.9)) and the normal space
\begin{equation}
N_{\sL(M)} BS( a_L) = H^1(\sL, \R).
\label{eq1.23}
\end{equation}
Now consider the space $hW\sL(M)_0$ of half-weighted Lagrangian cycles of zero volume (1.14) as a  bundle over $\sL(M)$. Then $hW\sL(M)_1$ is an affine bundle over it. Let us restrict these bundles to 
 the subspace $BS(L, a_L)$. We get the very one space 
\begin{equation}
hWBS( a_L)_t  \subset hW\sL(M)_t
\label{eq1.24}
\end{equation}
of half-weighted BS-Lagrangian cycles which is an  infinite-dimensional symplectic manifold with respect to the symplectic form $\Om$ (1.17). 

\begin{rmk} Recall that such symplectic form in the case $\sL = S^1$ is commonly used in connection with the Korteweg-de Vries equation (see \cite{DKN}). Such symplectic form has also arisen in the hamiltonian formulation of the motion of vortex patches for planar incompressible fluids. \end{rmk}

So we can repeat these constructions for Lagrangian submanifolds of  this infinite-dimension manifold (1.24). But first   we need to find a line bundle $L_B$ on $hWBS( a_L)_t$  with connection $A_W $ such that its curvature form is proportional to $\Om$ (1.17): 
\begin{equation}
F_A = 2 \pi \cdot K \cdot \Om , \quad K \in \Z.
\label{eq1.25}
\end{equation}

\subsection*{ From Lagrangian to Legendrian }

Consider the principal $U(1)$-bundle 
\begin{equation}
\pi \colon P \to M
\label{eq1.26}
\end{equation}
of our prequantization bundle $L$. The Hemitian connection $a_L$ is given by a 1-form $\al$. This form defines a contact structure on $P$ with the volume form
\begin{equation}
\frac{1}{2} \al \wedge d \al^n.
\label{eq1.27}
\end{equation}
If $\sL \in BS( a_L)$ then the restriction $(L, a_L) \vert_\sL$ admits  a covariant constant section
\begin{equation}
s_\sL  \colon \sL \to P \vert_\sL.
\label{eq1.28}
\end{equation}
Such a section is defined up to the natural $U(1)$-action on $P$ as  a principal
$U(1)$-bundle. 

Moreover, the submanifold 
\begin{equation}
s_\sL (\sL) \subset P
\label{eq1.29}
\end{equation}
is a {\it Legendrian} submanifold of the contact manifold $P$.
Such submanifolds of $P$ are called {\it Planckian}

 Let
\begin{equation}
\pi \colon  \La(P) \to BS( a_L) \subset \sL(M)
\label{eq1.30}
\end{equation}
be the space of all such Planckian  cycles over all BS-Lagrangian cycles where $\pi$ is the projection (1.26).

Obviously $\La(P)$ is a principal $U(1)$-bundle over $ BS( a_L)$.
\begin{prop} The tangent space of $\La(P)$ at a Planckian submanifold $s_\sL(\sL)$
\[
T\La(P)_{s_\sL(\sL)} = C^\infty (\sL).
\]
\end{prop}

\begin{dfn} The complex line bundle $L_B$ of the principal $U(1)$-bundle $\La(P)$ (1.30) over $ BS(L, a_L)$ is called the {\it Berry bundle}.
\end{dfn}

Consider the forgetful  map 
\begin{equation}
F \colon hWBS( a_L) \to  BS( a_L)
\label{eq1.31}
\end{equation}
as a bundle over $ BS( a_L)$ and lift it to the space of Planckian cycles $\La(P)$. We get the space $hW\La(P)$ of half-weighted Planckian cycles  and the forgetful map
\begin{equation}
F \colon hW\La(P)) \to  \La(P). 
\label{eq1.32}
\end{equation}
 
\begin{prop} This forgetful  map is the cotangent bundle of $\La(P)$:
\[
 hW\La(P) = T^*\La(P).
\]
\end{prop}
We get this statement immediately from (1.11) and (1.12).
\begin{cor} The space $ hW\La(P)$ is a symplectic infinite-dimensional manifold with respect to the  2-form 
\[
\Om_P = d \al
\]
where $\al$ is the standard action 1-form on a  cotangent bundle.
\end{cor}
Moreover the $U(1)$-action on $\La(P)$ defines the action of  $U(1)$ on the cotangent bundle $ T^*\La(P) =  hW\La(P)$. This action preserves the form $\Om_P$. 

\begin{prop}\begin{enumerate}
\item The function $v$ from (1.13) is a moment map for this $U(1)$-action. 
\item The symplectic quotient (the result of  the symplectic reduction procedure 
\[
v^{-1}(t) / U(1) = hWBS( a_L)_t
\]
is the space of half-weighted Bohr-Sommerfeld cycles of volume $t$. 
\end{enumerate}\end{prop}

\begin{cor} The space $ hWBS( a_L)_t$ is an infinite-dimensional symplectic manifold with the symplectic 2-form $\Om_t$ which is given by  the symplectic reduction procedure from the form $\Om_P$. The form $\Om_0$ is given by formula (1.17).
\end{cor}

Summarizing we get
\begin{enumerate} \item a 1-parameter family of phase spaces
\begin{equation}
( hWBS( a_L)_t, \Om_t)
\label{eq1.33}
\end{equation}
\item equiped with the Berry line bundle $L_B$ 
\item with principal $U(1)$-bundle
\[
\pi \colon F^* \La(P) \to  hWBS( a_L)_t
\]
lifted from (1.30).
\end{enumerate}

Now we need to find on $ hWBS( a_L)_t$ Hermitian connections on the Berry bundle $L_B$  with a curvature $\Om_t$.

It is easy to check the following 
\begin{prop} 
\begin{enumerate}
\item Every section
\[
S \colon  BS( a_L) \to  hWBS( a_L)_t
\]
defines an Hermitian connection $A_S$ on the Berry line bundle $L_B$ over $ BS( a_L)$;
\item the curvature form of such connection
\[
F_{A_S} = 2 \pi \cdot i \cdot \Om_t.
\]
\end{enumerate}
\end{prop}

A choice of such section $S$ is called a {\it  half-weighting rule}. Thus every half-weight rule defines the prequantised system
\begin{equation}
( BS( a_L)_t, S^* \Om_t, L_B, A_S )
\label{eq1.34}
\end{equation}
(which is ready to be quantized, see (1.19)). 

 The important half-weighting rule is a choice  an admissible (with $\om$) Riemannian metric and a metaplectic structure on $M$. Such metric defines the  admissible almost  complex structure. If this complex structure is integrable then $M$ becomes a {\it  Kahler manifold} and we can switch on strong  algebro-geometric methods.

\section{Complex structure }

To quantize a prequantized system $(M, \om, L, a_L)$ we need to fix upon some {\it polarization} of it. The complex polarization  is a choice of a complex
structure $I$ on $M$  such that $M_I$ is a K\"ahler manifold with K\"ahler form
 $\om$. Then the curvature form of the Hermitian connection $a$ is of type
$(1,1)$, hence  the line bundle $L$ is a
holomorphic line bundle on $M_I$.  Complex quantization provides the space
of  wave functions
\begin{equation}
\sH_I = H^0(L)
\label{eq2.1}
\end{equation}
-the space of holomorphic sections of the line bundle $L$.
 (see the survey \cite{K}).
\begin{rmk} Here we will consider the compact case only. In this case the wave function space (2.1) is finite-dimensional. But all our constructions are correct even in non compact case, for example when $M$ is the cotangent bundle of a configurations space where  the space of wave functions are  infinite-dimensional.
\end{rmk}

Returning to Proposition 1.3 we can see that the space $hW\La(P)$ admits  a complex structure $I_{P}$ as a cotangent bundle: the tangent space at an half-weighted Planckian $\La \subset P$
\begin{equation}
ThW\La(P)_{\La, hF} = C^\infty (\pi(\La)) \oplus  C^\infty (\pi(\La)) = C^\infty_\C (\pi(\La)).
\label{eq2.2}
\end{equation}
It is easy to check the following 
\begin{prop}
\begin{enumerate}
\item This complex structure $I_P$ of (2.2) is integrable.
\item This complex structure is compatible with the symplectic structure $\Om_P$ (see Corollary 1.1).
\end{enumerate}
\end{prop}
Thus our space $hW\La(P)$ is a Kahler infinite-dimensional manifold. Moreover this Kahler structure is invariant with respect to $U(1)$-action on $hW\La(P)$ induced by $U(1)$-action on $\La(P)$. Hence we can switch on the symplectic reduction with the moment map (1.13) (see Proposition 1.4).  
 
\begin{cor}
\begin{enumerate}
 \item  The symplectic quotient
\[
v^{-1}(t) = hWBS( a_L)_t
\]
admits a  quotient Kahler structure.
\item Every quadruple (1.33) of the  family of phase spaces is a Kahler manifold and admits a  complex polarization given by the symplectic quotient of the complex structure $I_P$ (2.2).
\item The Berry line bundle $L_B$ is holomorphic.
\item For every half-weighting rule 
\[
S \colon BS( a_L) \to hWBS( a_L)_t
\]
(see 1) of Proposition 1.5) the connection $A_S$ (1.34) is given by the holomorphic structure on the Berry line bundle $L_B$.
\item Every prequantized system of the family (1.34) admits the canonical complex (Kahlerien) polarization and can be quantized.
\end{enumerate}
\end{cor}

\begin{rmk} We would like accentuate that the  complex and Kahlerien structures  described here  don't depend on a choice of  Kahler structure $M_I$ on $M$.
\end{rmk}
  
Recall that $\pi(\La)$ is a BS-Lagrangian cycle in $M$. This cycle,  as any Lagrangian cycle,  {\it can't be contained by an algebraic divisor}. Thus the restriction map
\begin{equation}
res \colon  H^0(L) \hookrightarrow  C^\infty_\C (\pi(\La))
\label{eq2.3}
\end{equation} 
is an embedding. Therefore we have the trivial complex subbundle of the holomorphic  cotangent space of $hW\La(P)$:
\begin{equation}
0 \to \widetilde{ H^0(L)} \to T^*\La(P)
\label{eq2.4}
\end{equation}
or the trivial quotient
\begin{equation}
j \colon  T \La(P) \to \widetilde{ H^0(L)}^* \to 0.
\label{eq2.5}
\end{equation}
(Here $\widetilde{V}$ is the trivial vector bundle with a fiber $V$).
This epimorphism looks like a differential of a map to the vector space
$\sH^*=H^0(L)^*$. It is indeed.
Let us fix a metaplectic structure on $M_I$ (see for example section 1 of \cite{T1}). Then,  following Bortwick-Paul-Uribe \cite{BPU},
  we can construct the map 
\begin{equation}
BPU \colon hW\La(P) \to \sH^*
\label{eq2.6}
\end{equation}
with the differential (2.5). The  projectivization of this map is
\begin{equation}
\PP BPU \colon hWBS( a_L) \to \PP  \sH^*.
\label{eq2.7}
\end{equation}
 The construction is the  following:
\begin{enumerate}
\item our principal bundle $P$ (1.26) becomes the boundary of  the unit disc bundle $D \subset L^*$ 
\begin{equation}
\p D = P, \quad D \subset L
\label{eq2.8}
\end{equation}
where $D$  is a strictly pseudoconvex domain;
\item  there is the Szego orthogonal  projector 
\begin{equation}
\Pi \colon L^2_\C(P) \to H^0(L)
\label{eq2.9}
\end{equation}
 to the first component of  the Hardy space of boundary values of holomorphic   functions (in $D$) which are linear fiberwise;
\item for every Planckian submanifold $\La \in P$ in $C^\infty_\C(\La)$ there is the space of Legendrian distributions of order $m$ associated with $\La$ which is the Szego projection of the space of distributions  conormal  to $\La$     of order $m + \frac{dim M}{2}$ (see 2.1 of section 2 of \cite{BPU});
\item a half-form on $\La$ is identified with  the symbol of a Legendrian distribution of order $m$ (see 2.2 of section 2 of \cite{BPU}), thus at the symbolic level all Legendrian distributions look like delta-functions or their derivatives;
\item for a Legendrian submanifold  $\La$  equipped  with a half-form $hF$ we fix the Legendian distribution of order $\frac{1}{2}$ with symbol $hF$ which is the Szego projection of the delta function $\delta_\La$.  
\end{enumerate}
Summarizing
\begin{enumerate}
\item  for every Planckian submanifold  $\La \subset P$  equipped with a half-form $hF$ we have the vector
\begin{equation}
 BPU (\La, hF) = \Pi_{hF} (\delta_\La) \in H^0(L)
\label{eq2.10}
 \end{equation}
where $\Pi_{hF}$ is the Szego projection $\Pi$  to  the Hardy space  of the  distribution with the symbol $hF$ ;
\item for every Planckian $\La$ the image $\pi(\La)$ is a Bohr-Sommerfeld Lagrangian and every Planckian with this image is a translation of $\La$ by  $U(1)$-action on $P$, thus a pair 
$(\sL, hF)$ defines a point of the projectivization
\begin{equation}
 BPU (\sL, hF) = \PP ( \Pi_{hF} (\delta_\La)) \in \PP H^0(L)^*.
\label{eq2.11}
 \end{equation}
\end{enumerate}
Actually the Hermitian structure on $L$ defines the Hermitian structure on $ H^0(L)$. We can thus identify 
\[
 H^0(L) = \ov{ H^0(L)}^*.
\]
Our observations are the following

\begin{thm} \begin{enumerate}
\item The BPU-map (2.6) is holomorphic with respect to the complex structure from (2.2).
\item The epimorphism $j$ (2.5) is the differential 
\begin{equation}
j = d BPU.
\label{eq2.12}
 \end{equation}
\item The Berry bundle on $BS( a_L)$ is holomorphic and 
\begin{equation}
L_B = BPU^* (\Oh_{\PP H^0(L)^*} (1)).
\label{eq2.13}
 \end{equation}
\item Let $\rk H^0(L) = r_L$ and $S^{2r_L - 1}$ be the unit sphere in $ H^0(L)$. Then the principal $U(1)$-bundle 
\begin{equation}
F^* (\La(P)) = BPU^{-1}(S^{2r_L - 1})
\label{eq2.14}
 \end{equation}
and this equality is equivariant with respect to $U(1)$-action.
\item Let $\Om_{FS}$ be the Kahler form on $\PP H^0(L)$ given by the Fubini-Study metric induced by the identification $ H^0(L) = \ov{ H^0(L)}^*$. Then on
 $hWBS( a_L)_t$
\begin{equation}
\Om_t = BPU^* \Om_{FS}.
\label{eq2.15}
 \end{equation}
\item Let $A_{FS}$ be the connection on the Hopf bundle on $\PP H^0(L)^*$ induced by the Fubini-Study metric. Then for any half-weighting rule (see 1) of Proposition 1.5)
\begin{equation}
BPU^* A_{FS} \vert_{S(BS( a_L))}  
\label{eq2.16}
 \end{equation}
is a connection on the Berry bundle on $BS( a_L)$ with  the curvature form $\Om_t$.

\end{enumerate}
\end{thm}

 \section{Supercycles}

There is another way to give any Lagrangian cycle $\sL$ some  additional structure. This structure is a {\it  supercycle} structure or a {\it brane} structure with the support  $\sL$. 

\begin{dfn} A  pair $(\sL, a)$ where $a$ is a flat $U(1)$-connection of the trivial line bundle on $\sL$  is called a supercycle  ( or brane )  supported by  $\sL$.
\end{dfn} 
The trivial line bundle contains the trivial connection given by its covariant derivative 
\begin{equation}
\na_0 = d.
\label{eq3.1}
\end{equation} 
Then every connection $a$ on the trivial line bundle can be identify with 1-differential form $ i \cdot\al \in i \cdot  \Om(\sL)$. The flatness means that this form is closed.

Let the moduli space of supercycles be $S\sL(M)$. The the forgetful  map\begin{equation}
F \colon S\sL(M) \to \sL(M)
\label{eq3.2}
\end{equation} 
from this moduli space can be identified with the projection of the  tangent bundle. 
\begin{equation}
 S\sL(M) = T \sL(M).
\label{eq3.3}
\end{equation} 
Thus this space admits the natural complex structure: the tangent space\begin{equation}
T  S\sL(M)_{(\sL, a)} = Z^1(\sL) \oplus i \cdot Z^1(\sL) .
\label{eq3.4}
\end{equation} 
This almost complex structure is constant thus it is {\it integrable}.

So every infinitesiamal  deformation $(\al_1, i \cdot  \al_2)$ of a supercycle $(\sL, a)$ defines the  {\it  complex deformation} 
\begin{equation}
a + \al_1 + i \cdot \al_2
\label{eq3.5}
\end{equation}
of the flat connection $a$. The complex gauge class of this {\it complex}  connection defines a character of the fundamental group
\begin{equation}
\chi_{a + \al_1 + i \cdot \al_2} \colon \pi_1(\sL) \to \C^* \in H^1(\sL, \C) / H^1(\sL, \Z) = J^\C_\sL.
\label{eq3.6}
\end{equation}
Recall that there exists the map
\begin{equation}
r \colon  J^\C_\sL \to J_\sL , \quad r (\chi) = \frac{\chi}{\vert \chi \vert} 
\label{eq3.7}
\end{equation}

\begin{dfn} A deformation $(\al_1, i \cdot  \al_2)$ is called {\it  complex isodrastic} deformation of a supercycle $(\sL, a)$ if the character (3.7) is constant.
\end{dfn} 

The space  of isodrastic  infinitesimal deformations of a supercycle $(\sL, a)$  is
\begin{equation}
T_I S\sL(M)_{(\sL, a)} = B^1_\C = C^\infty_\C (\sL) /\C
\label{eq3.8}
\end{equation} 
(compair this formula with (1.19)).

Now every prequantization (1.19)-(1.20) equips every Lagrangian cycle $\sL$ with a supercycle structure 
\[
(\sL, \al_L \vert_\sL) \in S\sL(M)
\]
and we get the section
\begin{equation}
S_{a_L} \colon \sL(M) \to S\sL(M).
\label{eq3.9}
\end{equation}

Every supercycle $(\sL, a)$ defines the family of complex flat connections
\begin{equation}
a_L \vert_\sL + u \cdot a , \quad u \in \C.
\label{eq3.10}
\end{equation}
A complex gauge equivalence class of such  complex flat connections is
\begin{equation}
[a_L \vert_\sL + u \cdot a ] \in J^\C_\sL
\label{eq3.11}
\end{equation}

\begin{dfn} 
A Lagrangian supercycle $(\sL, a)$ is called complex Bohr-Sommerfeld ($BS_\C$ for short) if there exists $u \in \C$ such that
\begin{equation}
[a_L \vert_\sL + u \cdot a] = 0.
\label{eq3.12}
\end{equation}
\end{dfn}

Let 
\begin{equation}
BS_\C(a_L) \subset S\sL(M)
\label{eq3.13}
\end{equation}
be the space of all $BS_\C$-supercycles.
 In particular
\begin{equation}
(\sL, 0) \quad \text{ is } \quad  BS_\C \quad  \implies \quad \sL \in BS(a_L).
\label{eq3.14}
\end{equation}

On the other hand if $\sL \in BS(a_L)$ then
\begin{equation}
(\sL, a) \quad \text{ is } \quad  BS_\C \quad  \implies \quad a \in i \cdot B^1(\sL).
\label{eq3.15}
\end{equation}
So the space $BS_\C(a_L)$ contains the subspace
\begin{equation}
BS_\C(a_L)_0 = \{ (\sL, a) \quad \vert \quad \sL \in  BS(a_L), \quad a \in i \cdot B^1(\sL) \}. 
\label{eq3.16}
\end{equation}

The tangent space of $BS_\C(a_L)_0$ at $(\sL, a)$ is
\begin{equation}
(T BS_\C(a_L)_0)_{(\sL, a)} = B^1(\sL)_\C = C^\infty_\C (\sL) / \C.
\label{eq3.17}
\end{equation}

\begin{prop} If $u \neq u'$ and 
\begin{equation}
[a_L \vert_\sL + u \cdot a] = [a_L \vert_\sL + u' \cdot a] = 0 
\label{eq3.18}
\end{equation}
then 
\[
(\sL, a) \in  BS_\C(a_L)_0
\]
\end{prop}

On the space $ BS_\C(a_L)$ there thus exists the function
\begin{equation}
u \colon  BS_\C(a_L) \to \C
\label{eq3.19}
\end{equation}
with   level complex hypersurfaces of this function
\begin{equation}
 BS_\C(a_L)_u  \subset  BS_\C(a_L).
\label{eq3.20}
\end{equation}
\begin{rmk} Recall that the definition of the zero-level is given in (1.16).
\end{rmk}

\begin{prop} 
\begin{enumerate}
\item The tangent space of $BS_\C(a_L)$ at a supercycle $(\sL, a)$ is 
\begin{equation}
T BS_\C(L, a_L) = C^\infty_\C(\sL).
\label{eq3.21}
\end{equation}

\item The normal bundle
\begin{equation}
(N BS_\C(L, a_L)_u)_{(\sL, a)} = H^1(\sL,\C).
\label{eq3.22}
\end{equation}
\end{enumerate}
\end{prop}
We can see that this normal bundle doesn't depend on $a$ thus to be $BS_\C$ supercycle is a property of a Lagrangian cycle itself.

\begin{cor} The space $BS_\C( a_L)$  possess  a  complex structure and the 
function $u$ from (3.19) on it is holomorphic. \end{cor}

\subsection*{$\C^*$-lifting}
Consider the principal $\C^*$-bundle  
\begin{equation}
\pi \colon P_\C \to M
\label{eq3.23}
\end{equation}
of the prequantization line bundle $L$. Our connection $a_L$ is given by the complex  1-form $\al$ and the complex 2-form $d \al$ defines a complex  symplectic structure on $P_\C$:
\begin{equation}
 (P_\C, d  \al ).
\label{eq3.24}
\end{equation}

If $(\sL, a) \in BS_\C( a_L)$ then there exists a  covariant constant sections
\begin{equation}
s_{(\sL, a)}  \colon \sL \to P_\C  \vert_\sL.
\label{eq3.25}
\end{equation}
Such a section is defined up to the natural $\C^*$-action on $P_\C$ as on a principal $\C^*$-bundle. 

Moreover, the submanifold 
\begin{equation}
s_{(\sL, a)} (\sL) \subset P_\C
\label{eq3.26}
\end{equation}
is an isotropic  submanifold with respect to the 2-form from (3.24).
Such submanifolds of $P_\C$ are called {\it  Planckian} again.

 Let
\begin{equation}
\pi \colon  \La(P_\C) \to BS_\C( a_L) 
\label{eq3.27}
\end{equation}
be the space of all such complex  Planckian  cycles over  all $BS_\C$-Lagrangian cycles where $\pi$ is the projection (3.23).

Obviously (3.27) is a principal $\C^*$-bundle over $BS_\C( a_L)$.
\begin{prop} The tangent space of $\La(P_\C)$ at a Planckian submanifold $s_\sL(\sL)$
\[
T\La(P_\C)_{s_\sL(\sL)} = C^\infty_\C (\sL).
\]
\end{prop}

\begin{dfn} The complex line bundle $L_B$ of the principal $\C^*$-bundle
 $\La(P_\C)$ (3.27) over $BS_\C( a_L)$ is called the {\it Berry bundle}.
\end{dfn}

Consider the forgetful  map 
\begin{equation}
F \colon BS_\C( a_L) \to  \sL(M)
\label{eq3.28}
\end{equation}
as a bundle over $F(BS_\C( a_L)) $ and lift it to the space of complex  Planckian cycles $\La(P_\C)$. We get the space $S\La(P)$ of complex  Planckian supercycles  and the forgetful  map
\begin{equation}
F \colon S\La(P_\C)) \to  \La(P_\C). 
\label{eq3.29}
\end{equation}
 
\begin{prop} This forgetful  map is the tangent bundle of $\La(P_\C)$:
\[
 S\La(P_\C) = T\La(P_\C).
\]
\end{prop}
It is just an interpretation of previous constructions.

Moreover the $\C^*$-action on $\La(P_\C)$ defines a  holomorphic action of  $\C^*$ on the tangent bundle $ T\La(P_\C) $. 

 \begin{prop}\begin{enumerate}
\item The function $u$  (3.19)  is a moment map for this $\C^*$-action. 
\item The complex  quotient (the result of  the reduction procedure 
\[
 BS_\C( a_L)_u
\]
is the space of complex Bohr-Sommerfeld cycles of level $u$. 
\end{enumerate}\end{prop}

\begin{cor} The space $BS_\C(L, a_L)_u$ is an infinite-dimensional complex  manifold. 
\end{cor}

Summarizing we get
\begin{enumerate} \item a  1-complex parameter family of complex  spaces 
\begin{equation}
( BS_\C( a_L)_u )
\label{eq3.30}
\end{equation}
\item equipped with the holomorphic  Berry line bundle $L_B$ 
\item with principal $\C^*$-bundle
\[
\pi \colon F^* \La(P_\C) \to  BS_\C( a_L)_u.
\]
\end{enumerate}

Again  let us  fix  some   complex
structure $I$ on $M$  such that $M_I$ is a K\"ahler manifold with K\"ahler form
 $\om$. Then   the line bundle $L$ is a
holomorphic line bundle on $M_I$ with  the space
of  holomorphic sections $ H^0(L)$ (2.1).
  
Recall that $\pi(\La)$ is a BS-Lagrangian cycle in $M$.  Thus the restriction map
\begin{equation}
res \colon  H^0(L) \hookrightarrow  C^\infty_\C (\pi(\La))
\label{eq3.31}
\end{equation} 
is an embedding. Therefore we have the trivial complex subbundle of the holomorphic  tangent space of $\La(P_\C)$:
\begin{equation}
0 \to \widetilde{ H^0(L)} \to T \La(P_\C).
\label{eq3.32}
\end{equation}

This monomorphism defines the space of commuting holomorphic automorphisms: for every $s \in H^0(L)$ and $ (\sL, a) \in BS_\C(a_L)_u $
\begin{equation}
t_s  \colon \La(P_\C) \to \colon \La(P_\C) 
\label{eq3.33}
\end{equation}
commute with the $\C^*$-action. Thus we have the family of holomorphic automorphisms
\begin{equation}
t_s  \colon BS_\C(a_L) \to  BS_\C(a_L)
\label{eq3.34}
\end{equation}
preserving the holomorphic Berry line bundle $L_B$.

\begin{cor} The space $H^0(L)$ acts on the space $H^0(L_B)$ of holomorphic sections of the Berry line bundle. This action is homogeneous with respect to the $\C^*$-action.
\end{cor}.  

\subsection*{ Complex structure }

Now suppose that our  $M$ is equipped with a Riemannian metric $g$ compatible with the symplectic structure $\om$.  For us the  important case is when  our symplectic manifold $(M, \om)$ is equipped with a compatible complex structure  $I$   which is integrable.  In this case $(M_I, \om)$ is a  Kahler manifold.  Then the space $\sA^0 = Z^1(\sL)$ of flat connections  on a Lagrangian submanifold $\sL$ has  the finite dimensional subspace
\begin{equation}
H^1(\sL, \R)_h = \{ i \cdot \al \vert d \al = 0,\quad d * \al = 0 \} \in \sA^0_0 
\label{eq3.35}
 \end{equation}
of harmonic forms with respect $g$. So the space $S\sL\sM$ of all  Lagrangian super cycles contains the subspace
\begin{equation}
S\sL(M_I)^h  \subset S\sL(M)
\label{eq3.36}
 \end{equation}
of harmonic super cycles (harmonic with respect to the Riemannian metric induced by the compatible complex structure $I$).

\begin{rmk} It is quite convenient to code a compatible Riemanian metric on a symplectic manifold by the corresponding almost complex structure. Usually we are considering Kahler manifolds where this  almost complex structure $I$ is integrable and the  metric is Kahlerian.
\end{rmk}

 The gauge group 
\[
\sG = Map (M, U(1))
\]  actes on our space $ S\sL(M) $ of all Lagrangian supercycles (on the second component of every pair $(\sL, a)$). The orbit  space 
\begin{equation}
QS = S\sL\sM / \sG
\label{eq3.37}
 \end{equation}
is fibered over the space of Lagrangian cycles 
\begin{equation}
F_Q  \colon QS  \to \sL(M)
\label{eq3.38}
 \end{equation}
with finite dimensional fibers. Such fiber  over a smooth Lagrangian submanifold $\sL$ is
\begin{equation}
F_Q^{-1} (\sL) =  H^1(\sL, \R) / H^1(\sL, \Z) = J_\sL.
\label{eq3.39}
 \end{equation}
is  the jacobian of $\sL$ (see (1.6). We can say that the fibration $F_Q$ is the {\it jacobian fibration of the universal Lagrangian cycle}.

The factorization map (3.40) restricted   to the subspace  (3.39) of harmonic connections  gives the surjection
\begin{equation}
U \colon S\sL(M_I)^h \to QS.
\label{eq3.40}
 \end{equation}
 Over a smooth $\sL$ this map  is the universal cover
\begin{equation}
U \colon H^1(\sL, \R) \to J_\sL.
\label{eq3.41}
 \end{equation}
The intersection 
\begin{equation}
BS_\C(a_L) \cap S\sL(M_I)^h = BS_\C(a_L)^h
\label{eq3.42}
 \end{equation}
is the space of harmonic $BS_\C$-cycles. In particular for every value of the function $u$ (3.19) we have the space 
\begin{equation}
BS_\C(a_L)_u \cap S\sL(M_I)^h = BS_\C(a_L)^h_u
\label{eq3.43}
 \end{equation}
of harmonic supercycles of level $u$. In particular it is easy to see that

\begin{prop}
\[
BS_\C(a_L)^h_0 = BS(a_L) 
\]
\end{prop}

Obviously  we can map these spaces to $QS$ (see (3.43))
\begin{equation}
U \colon BS_\C(a_L)^h_u \to QBS_\C(a_L)^h_u \subset QS.
\label{eq3.44}
 \end{equation}
 Over a smooth $\sL$ this map  is the universal cover (3.44).

Now we suppose that $M_I$ admits a metaplectic structure (see section 1 of \cite{T1}). Then every smooth Lagrangian submanifold admits a metalinear structure. Any  Kahler metric $g$ on $M$ defines a half-weighting rule 
\begin{equation}
S_g \colon \sL(M) \to hW\sL(M)
\label{eq3.45}
 \end{equation}
(see Proposition 1.5) which equips every Lagrangian cycle $\sL$ with a half-form $hF$ such that 
\begin{equation}
hF^2 = \frac{V_g}{\int_\sL V_g}
\label{eq3.46}
 \end{equation}
where $V_g$ is the Riemannian volume form on $\sL$.

As any half-weighting rule the section (3.45) defines the identification $ T\La(P_\C = T^*\La(P_\C$ and defines Kahler structures on every complex quotient $BS_\C(a_L)_u$ of Proposition 3.5.

Moreover the BPU-construction (2.6)-(2.7) is working. Indeed to switch  it on  we have to define a Legendrian submanifold of $P$ (2.8) over a Lagrangian submanifold $\sL \subset M$. But for $\sL \in F(BS_\C(a_L)$ we have a lifting  
\begin{equation}
s_\sL \colon \sL \to  P_\C.
\label{eq3.47}
 \end{equation}
Then the projection 
\begin{equation}
\frac{s_\sL}{\vert s_\sL \vert} \colon \sL \to P
\label{eq3.48}
 \end{equation}
gives us what we need (see (3.7)). 

Now using the BPU-construction we get the map
\begin{equation}
BPU_\C  \colon \La(P_\C) \to H^0(L)^*
\label{eq3.49}
 \end{equation}
and the projectivization of it
\begin{equation}
\PP BPU_\C  \colon BS_\C(a_L)  \to \PP  H^0(L)^*.
\label{eq3.50}
 \end{equation}

Again we can check that

\begin{thm}
\begin{enumerate}
\item The map $BPU_\C$ (3.49) is holomorphic with respect to the complex structure of Proposition 3.5;

\item the epimorphism $res^*$ (see (3.31)) is the differential of $BPU_\C$; 
\item the Berry bundle
\begin{equation}
L_B = BPU_\C^*(\Oh_{\PP H^0(L)^*} (1));
\label{eq3.51}
 \end{equation}
\item and so on like in Theorem 2.1.
\end{enumerate}
\end{thm}

There are two strategies  to use this new complex parameter $u$:
\begin{enumerate}
\item consider the situation when $ u = \frac{1}{N}$ is a real rational number 
\[
  \quad \text{ and let  } \quad N \to \infty;
\]
\item use $u$ as a complex parameter to create get a geometric object- a complex algebraic curve - a  ``spectral curve''.\end{enumerate}

In the first case  equipping Lagrangian cycles with half-forms  and applying the direct BPU-method we get  special configurations of states in wave functions spaces $H^0(L^k)$. This gives an integer structure on vector spaces (2.1).

In this paper we are using the second strategy applying these constructions  to  examples considered in \cite{T1} and \cite{T2}.

\section{$u$-curves of  real  polarizations }

 Recall that we do not need a complex structure on $M$  to define   a  real polarization of  a quadruple 
 \begin{equation} 
( M, \om, L, a_L )
\label{eq4.1}
 \end{equation}
with the condition (1.19) on  the curvature of an Hermitian connection $a_L$.
Then a real polarization  of this quadruple 
 is a fibration
 \begin{equation}
 \pi\colon M \to B,
\label{eq4.2}
 \end{equation}
such that $ \om\rest{\pi^{-1}(b)}=0$ for every point $b\in B$
and for generic $b$ the fibre $\pi^{-1}(b)$ is a smooth Lagrangian.
If we consider the pair $(M,\om)$ as a phase space of a mechanical
system ( see for example \cite{A}) it admits a real polarization if and only if it is completely  integrable. In the compact case a generic fiber is a $n$-torus
$T^n$, $ (2n = dim_{\R} M)$ and $dim B = n$. As every family of Lagrangian cycles the family (4.2) defines the Kodaira-Spencer homomorphism
 \begin{equation}
 KS_b  \colon TB_b \to H^1(\pi^{-1}(b), \R)
\label{eq4.3}
 \end{equation}
at every smooth fiber $\pi^{-1}(b)$. It is a well know fact that if $M$ is compact then every smooth fiber 
\[
\pi^{-1}(b) = T^{\frac{dim M}{2}}
\]
is $n$-torus where $dim M = 2n$.

 We would like to add the following condition to the definition of a real polarization: for every smooth fiber $\pi^{-1}(b) = T^n$ the Kodaira-Spencer map (4.3)
 \begin{equation}
 KS_b  \colon TB_b \to H^1(T^n, \R)
\label{eq4.4}
 \end{equation}
is an isomorphism. 

Under this condition we can restrict the fibration $F_Q$ (3.38) to $B$. We get the finite-dimensional family
\begin{equation}
F_Q \colon QS \to B 
\label{eq4.5}
 \end{equation}
of Jacobians of fibers of the projection $\pi$ (4.2).

Now let us switch to  a complex structure $M_I$ on $M$ and the corresponding  Kahler metric $g$. For every level $u \in \C$ consider the intersection
\begin{equation}
SB \cap BS_\C(a_L)_u  
\label{eq4.6}
 \end{equation}
and the projection of this set to $QS$ over $B$ (4.5). We have to get a finite set of classes of supercycles:
\begin{equation}
F_Q ( SB \cap BS_\C(a_L)_u ) = \{ (\sL_1, [a_1])_u, ... , (\sL_N, [a_N])_u \} 
\label{eq4.7}
 \end{equation}
where $[a_i] \in J_{\sL_i}$ (see (3.22)).

 When $u$ swept out $\C  $ we get a affine algebraic curve
\begin{equation}
 \Si  = \{ (\sL_i, [a_i])_u \} \subset  BS_\C(a_{L^k})
\label{eq4.8}
 \end{equation}

This curve  can be compactified and normalized to a compact algebraic curve
\begin{equation}
\Si^k_{\pi, I}  
\label{eq4.9}
 \end{equation}

\begin{dfn} This curve is called a  $u$-curve of the  mechanical system
(4.1) with respect to a real polarization $\pi$ (4.2) and a complex polarization $I$.
\end{dfn}

There are two possibilities  for the  behavior of  such a curve with respect to the complex parameter $u$:
\begin{enumerate}
\item either the natural projection 
\begin{equation}
u \to \Si^k_{\pi, I}  
\label{eq4.10}
 \end{equation}
is the universal cover. Then $\Si^k_{\pi, I} $ is an elliptic curve;

\item or there is the pencil
\begin{equation}
\phi \colon  \Si^k_{\pi, I} \to \PP^1_u  
\label{eq4.11}
 \end{equation}
such that 
\begin{equation}
\phi^{-1}(u) = \{ (\sL_i, [a_i])_u \}   
\label{eq4.12}
 \end{equation}
is the effective divisor (4.7).
\end{enumerate}

In the second case the atributes of this curve are
\begin{enumerate}
\item the   pencil 
\begin{equation}
\vert \xi \vert = \PP^1_u
\label{eq4.13}
 \end{equation}
(may be with base points) which gives the finite cover (4.9).

\item \[ 
 \phi^{-1}(0) = B \cap BS^k(M, L)
\]
is the set of original Bohr-Sommerfeld fibers (see Example 1 of section 1 from \cite{T1});

\item the degree of the surjection (2.4)
\begin{equation}
deg \phi  = \# ( B \cap BS^k(M, L))
\label{eq4.14}
 \end{equation}
is the number of Bohr-Sommerfeld fibers of level $k$ of the fibration $\pi$ (4.2)
\item the map 
\begin{equation}
i_{\pi, I}  \colon \Si^k_\pi \hookrightarrow F_Q^{-1}(B)
\label{eq4.15}
 \end{equation}
which for a point $x \in \Si^k_{\pi, I} $ over $u$ gives the $u$-Bohr-Sommerfeld Lagrangian super  fiber of our polarization;

\item the line bundle $L_{\pi, I}$ on $\Si^k_{\pi, I} $ the fiber of which over a point
 $(\sL, [a]) \in \Si^k_{\pi, I} $ over $u$ is the line 
\[
\pi^{-1} (\sL) \in \La(P_\C)
\]
(see (3.23)).That is this line bundle 
\begin{equation}
L_{\pi, I} = res L_B
\label{eq4.16}
 \end{equation}
is the restriction of the Berry line bundle (see Definition 3.4).

\end{enumerate}

Summarizing let us fix a real polarization $\pi$ (4.2) of $(M, \om)$ from the quadruple (4.1). Then sending a complex structure $I$ to the curve $\Si^k_{\pi, I}$ (4.8) we get the map
\begin{equation}
m \colon \sM \to M_g
\label{eq4.17}
   \end{equation}
of the moduli space of polarized complex structures on $M$  to the moduli space of curves of genus $g$. 

Not many explicit examples can be shown here. Despite this lack 
 we consider   applications  of these  constructions to the subject and examples of paper \cite{T1}.

\subsection*{Elliptic curve}

 The classical theory of theta functions is the first beautiful subject for application  of our construction. Let us start with dimension one  case. 

 Let $E$ be an elliptic curve  with zero element $o\in E$ and with flat metric
 $g$. Then the tangent bundle $TA$ has the standard constant Hermitian
 structure (that is, the Euclidean metric, symplectic form and complex
 structure $I$). The K\"ahler form $\om$ gives a polarization of degree 1.
 If we switch to  a complex structure we get a phase space of a classical
mechanical system 
\begin{equation}
 (A = T^2, \om, L = \Oh_E(o), a )
\label{eq4.18}
 \end{equation}
just like (4.1). Let us fix a smooth Lagrangian decomposition of $E$
  \begin{equation}
 E= T^{2}=S^1_+\times S^1_-. 
\label{eq4.19}
 \end{equation}
 Circles of  both families  are
 Lagrangian with respect to $\om$.  This  decomposition 
 induces the Lagrangian  decomposition $ H^1(E,\Z)=\Z_+\times\Z_-,$  and the  Lagrangian  decomposition
$ E_k=(S^1_+)_k \times (S^1_-)_k $ of the group of points of any  order $k$ i. e. the collection of theta structures of any level $k$ (see \cite{Mum}).
Thus  this  Lagrangian  decomposition  defines the collection of
 decompositions  of the
spaces of
wave functions
 \begin{equation}
 \sH_I^k=H^0(E, L^k)=\bigoplus_{w\in (\Z)_k^-}\C\cdot\theta_w,
 \quad\text{with}\quad \rk\sH_I^k=k,
 \label{eq4.20}
 \end{equation}
where $\theta_w$ is the theta function with  characteristic $w$ (see
\cite{Mum}).

On the other hand the  direct product
(\ref{eq4.19}) gives us a real polarization
 \begin{equation}
\pi\colon E\to S^1_-=B.
 \label{eq4.21}
 \end{equation}
Remark that in this case the action coordinates  (see for example \cite{A}) are just the flat coordinates  on $S^1_-=B$, and under this identification
 \begin{equation}
 B \cap BS^k = (S^1_-)_k
 \label{eq4.22}
 \end{equation}
is the subgroup of points of order $k$.

Thus for our  $u$-curve $\Si^k_{\pi, I}$ (4.5) in this case 
 \begin{equation}
\phi^{-1}(0) = (S^1_-)_k .
 \label{eq4.23}
 \end{equation}

Moreover, the embedding $i_{\pi, I}$ (4.8) sends  the collection of points
\[
\phi^{-1}(0) \subset \Si^k_{\pi, I}
 \]
to $B = S^1_-$ and precisely
\begin{equation}
 i_{\pi, I} (\phi^{-1}(0)) =  (S^1_-)_k \in  S^1_- .
 \label{eq4.24}
 \end{equation}

Now 
\begin{equation}
\pi' \colon E' = SQ \to S^1_-=B
 \label{eq4.25}
 \end{equation}
is the mirror partner of $E$ (see the diagram (1.49) from \cite{T2}).

Summarizing we have

\begin{prop}
\begin{enumerate}
\item The $u$-curve $ \Si^k_{\pi, I}$ of the system (4.18) is an elliptic curve;
\item the map (4.15) is the isogeny of order $k$
\begin{equation}
i_{\pi, I} \colon \Si^k_{\pi, I} \to E'
 \label{eq4.26}
 \end{equation}
induced by the isogeny 
\[
 \mu_k  \colon  S^1_- \to  S^1_-.
\]
We have thus the first case (4.10);

\item the map $m$ (4.17) 
\[
m \colon M_1 \to M_1
\]
is the isogeny map along the period $S^1_-$;

\item the Berry bundle (4.16) is 
\begin{equation}
L_{\pi, I} =  i_{\pi, I}^* \Oh_{E'}(o) 
 \label{eq4.27}
 \end{equation}
the line bundle of degree $k$.
\end{enumerate}
\end{prop}

Finally the canonical class of $E$ is zero. Thus $E$ admits a metaplectic structure and every circle admits a metalinear structure. Thus  a Kahler metric $g$ on $E$ defines a half-weighting rule (3.45) and $BPU_\C$ map
\begin{equation}
BPU_\C \colon \Si^k_{\pi, I} \to \PP \sH^k_I 
 \label{eq4.28}
 \end{equation}
(see (4.20)).
This is the standard embedding by the complete linear system of degree $k$ 
to the space with the basis (4.20).

\section{Non-abelian theory of theta functions}

Now  we would like to  apply these constructions to the theory of  non commutative
theta-functions (see \cite{T1}). To come to our standard setup 
consider the (6g-6)-manifold
 \begin{equation}
 R_g = \Hom (\pi_1(\Si_g),\SU(2)) /\PU(2)
 \label{eq5.1}
 \end{equation}
-- the space of classes of $\SU(2)$-representations of the fundamental
group of the  Riemann surface $\Si_g$ of genus $g$. This space can be included to the quadruple
 \begin{equation}
 (R_g,\Om, L, A_{\CS})
 \label{eq5.2}
 \end{equation}
where the symplectic  form $\Om$ can be constructed directly following W. Goldman or  using symplectic reduction
arguments as in \cite{RSW} where the determinant line bundle $L$ with the Chern-Simons connection were constructed also.  By the construction, the curvature form of this connection is
\begin{equation}
 F_{A_{\CS}}=i\cdot\Om.
  \label{eq5.3}
 \end{equation}
So we are in the very one situation of the Geometric Quantization.

To switch on a complex structure on $R_g$ we have to give a conformal structure $I$ to our  Riemann
surface $\Si$ of genus $g$. We get a complex
structure on the space of classes of representations $R_g$ such that
\begin{equation}
 R_{\Si}=R_g=\sM\semis
 \label{eq5.4}
 \end{equation}
 is the moduli space of semistable holomorphic
vector bundles on $\Si$. Then the form $F_{A_{\CS}}$ (5.3) is a $(1,1)$-form and
the line bundle $L$ admits the unique holomorphic structure compatible with
the Hermitian connection $A_{\CS}$. Thus we get the system of wave functions spaces
 \begin{equation}
 \sH_I^k = H^0 (R_I, L^k).
 \label{eq5.5}
 \end{equation}
 Ranks of these spaces are
given by the Verlinde formula (See  \cite{B} ):
 \begin{equation}
 rk \quad \sH_I^k = \frac{(k+2)^{g-1}}{2^{g-1}} \sum_{n=1}^{k+1}
 \frac{1}{(\sin(\frac{n\pi}{k+2}))^{2g-2}}\,.
 \label{eq5.6}
 \end{equation}

A real polarization of $R_g, \Om$ is defined by a trinion decomposition of our Riemann surface $\Si$. Such a decomposition is  given by the  choice of a
maximal collection of disjoint, noncontractible, pairwise nonisotopic
smooth circles on $\Si$. It is easy to see
 that any such system contains $3g-3$ simple closed circles
 \begin{equation}
C_1, \dots, C_{3g-3}\subset \Si_g,
 \label{eq5.7}
 \end{equation}
and the complement is the union
of $2g-2$ trinions $P_j$.

 The invariant of such a decomposition is given by its {\em $3$-valent
dual graph} $\Ga(\{C_i\})$,
associating a vertex to each trinion $P_i$, and an edge linking $P_i$ and
$P_j$ to a circle $C_l$  such that
 \[
 C_l\,\subset\,\p P_i \cap\p P_j.
 \]

 The isotopy type of the system (5.7) defines the map
 \begin{equation} \pi_{\{C_i\}} \colon R_g\to \R^{3g-3}
 \label{eq5.8}
 \end{equation}
to the real euclidian space with fixed coordinates $(c_1, \dots, c_{3g-3})$ such that
 \begin{equation}
 c_i (\pi_{\{C_i\}} (\rho))=\frac{1}{\pi}\,\cos^{-1}(\frac{1}{2}\tr
\rho([C_i]))\in [-1, 1].
 \label{eq5.9}
 \end{equation}
where $ \{ C_i \} = E(\Ga)$ - the set of edges of our graph. Then (see \cite{JW}) the map $\pi_{\{C_i\}}$ is a real polarization of the system
$(R_g,k\cdot\om,L^k,k\cdot A_{\CS})$and  coordinates $c_i$ are  action coordinates for this Hamiltonian system.

Moreover  the image of $R_g$ under $\pi_{\{C_i\}}$
 is a convex polyhedron 
\begin{equation}
 \Delta_{\{C_i\}}\subset [0,
 1]^{3g-3}. 
 \label{eq5.10} \end{equation}
Beside of this polyhedron our  space contains 
the integer sublattice $\Z^{3g-3}\subset \R^{3g-3}$, and we can consider
the " action" torus:
  \begin{equation} T^A=\R^{3g-3} / \Z^{3g-3}
\label{eq5.11} \end{equation}
containing the   topological complex
 \begin{equation}
 \ov{\Delta_{\{C_i\}}}
 \label{eq5.12}
 \end{equation}
which is the image of $ \Delta_{\{C_i\}}$ (5.10) in this action torus.
Then we have the following description of Bohr-Sommerfeld fibers of level $k$: first of all the intersection 
 \begin{equation}
BS^k \cap  \ov{\Delta_{\{C_i\}}} =  \ov{\Delta_{\{C_i\}}} \cap (T^A)_k 
 \label{eq5.13}
 \end{equation}
where $(T^A)_k$ is the subgroup of points of order $k$ of the action-torus.

To describe the inquired subset  we consider our  3-valent graph $\Ga$ with the set  $V(\Ga)$ of
vertexes ( $\vert V(\Ga) \vert = 2g-2$) and the set $E(\Ga)$ of edges
($\vert E(\Ga) \vert = 3g-3$) and the set $W^k_g$ of all  functions
 \begin{equation}
 w\colon E(\Ga) \to \{0, \frac{1}{k}, ... , \frac{k-1}{k} \}
 \label{eq5.14}
 \end{equation}
on the collection of edges of the 3-valent graph $\Ga $ to the
collection of $\frac{1}{k}$ integers. Then
\[
W^k_g = (T^A)_k.
\]
This set contains the subset $W^k_g(\Ga)$ of functions subjecting to conditions:  for any three edges
$C_l,
C_m, C_n$ meeting at a vertex $P_i$
 \begin{enumerate}
 \item $w(C_l) + w(C_m) + w(C_n)\in \frac{2}{k}\cdot \Z$;

 \item $w(C_l) + w(C_m) + w(C_n) \le 2$;

 \item for any ordering of the triple $C_l, C_m, C_n$,
 \begin{equation}
 \vert w(C_l)-w(C_m) \vert \le w (C_n) \le w(C_l) + w(C_m);
 \label{eq5.15}
 \end{equation}
\item if an edge $C_i$ separates the graph then
\[
w(C_i) \in \frac{2}{k} \cdot \Z
\]
\end{enumerate}
Such  function $w$ is called an   {\it admissible integer weight of level
$k$}  on the graph $\Ga$ ( see \cite{JW}). Now the statement of the proposition 
follows  from results of \cite{JW}.

 Let us apply our constructions from the previous sections to this situation. We get the  $u$-curve 
 \begin{equation}
\Si^k_{\pi, I}  \subset BS_\C(a_{L^k})
\label{eq5.16}
 \end{equation}
(see Definition 4.1). It is easy to see that here we have the case 2 of the alternative after Definition 4.1. Thus this curve has the full collection of atributes (4.11) - (4.16). Firstly we need two of them:

the embedding (4.15)
\begin{equation}
i_{\pi, I} \colon    \Si^k_{\pi, I}  \to F^{-1}_Q (\Delta_{\{C_i\}})
\label{eq5.17}
 \end{equation}
and the complex BPU-map (3.50)
 \begin{equation}
\PP BPU_\C \colon \Si^k_{\pi, I} \to \PP H^0(L^k)^*.
\label{eq5.18}
 \end{equation} 
Recall that the space $ F^{-1}_Q (\Delta_{\{C_i\}})$ (5.16) contains the Geometric Fourier Transformation of $L^k$
 \begin{equation}
GFT (L^k) \subset  F^{-1}_Q (\Delta_{\{C_i\}}).
 \label{eq5.19}
 \end{equation}
The intersection
\begin{equation}
GFT (L^k) \cap i_{\pi, I}(\Si^k_{\pi, I}) = BS(a_{L^k})
 \label{eq5.20}
 \end{equation}
-the collection of  $BS^k$-fibers. Moreover
\begin{equation}
BPU_\C (GFT (L^k) \cap i_{\pi, I}(\Si^k_{\pi, I})) \subset \PP H^0(L^k)^* 
 \label{eq5.21}
 \end{equation}
is the Bohr-Sommerfeld basis in $ H^0(L^k)^* $ (see \cite{T1}).

The combination of map (5.16) and the forgetful  map $F_Q$ defines the map
\begin{equation}
F \circ i_{\pi, I} \colon  \Si^k_{\pi, I} \to  \ov{\Delta_{\{C_i\}}} \subset T^A.
 \label{eq5.22}
 \end{equation}
This map induces the Jacobians  homomorphism
\begin{equation}
 j_{\pi, I}  \colon J_{\Si^k_\pi} \to   T^A .
 \label{eq5.23}
 \end{equation}
\begin{prop} This homorphism of real tori is surjective
\[
 J_{\Si^k_\pi} \to   T^A  \to 0
\]
\end{prop}
It is follows immediately  from (5.12)-(5.13).

\begin{prop} There exists the direct decomposition 
\[
 J_{\Si^k_\pi} = T_0 \times T^{3g-3}_+ \times T^A
\]
on real tori such that the surjection $j_\pi$ is the composition
\[
 J_{\Si^k_\pi} \to  T^{3g-3}_+ \times T^A \to T^A
\]
\end{prop}

\begin{cor} \begin{enumerate}
\item This decomposition gives the  decomposition of the subgroup of order $k$ points 
\[
(J_{\Si^k_{\pi, I}})_k = (T_0)_k \times (T^{3g-3}_+)_k \times (T^A)_k
\]
  \item and the embedding
\[
(T^A)_k \hookrightarrow (0, 0, (T^A)_k) \subset (J_{\Si^k_{\pi, I}})_k 
\]
\item and hence the embedding
\[
W^k_g (\Ga) \hookrightarrow (J_{\Si^k_{\pi, I}})_k.
\]
\end{enumerate} \end{cor}

Now our main geometrical observation is as  follows
\begin{prop}
The collection of non-abelian theta functions with characteristics ( elements of Bohr-Sommerfeld basis) of level $k$  corresponds  to the subcollection of abelian theta functions with characteristics of the $u$-curve $\Si^k_{\pi, I}$.
\end{prop}

\begin{prop} There exists a standard embedding of the $u$-curve in  the jacobian\[
    \Si^k_\pi \hookrightarrow J_{\Si^k_{\pi, I}}
\]
such that the intersection
\[ 
 \Si^k_{\pi, I} \cap W^k_g = W^k_g (\Ga) =  \ov{\Delta_{\{C_i\}}} \cap BS(a_{L^k}).
\]
\end{prop}

We have to stop here the development of a ``theory'' that isn't accompanied by a large collection of examples yet.

\subsection*{Acknowledgments}

I would like to express my gratitude to the University of  Oslo  and personally to Kristian Ranestad  and Geir Ellingsrud  for
support and hospitality.

\bigskip
\noindent
Andrei Tyurin, Algebra Section, Steklov Math Institute,\\
Ul.\ Gubkina 8, Moscow, GSP--1, 117966, Russia \\
e-mail: Tyurin@tyurin.mian.su {\em or} Tyurin@Maths.Warwick.Ac.UK\\
{\em or} Tyurin@mpim-bonn.mpg.de
 \end{document}